\documentclass[11pt]{article}
\usepackage{graphicx, color, xcolor}

\usepackage{amsmath,amsfonts,amssymb}
\usepackage{amsthm}

\def\forall{\hbox{for all}~}
\def\L{\mathbf{L}}

\def\ve{\varepsilon}

\def\D{{\cal D}}

\def\R{{\mathbb R}}

\def\implies{\Longrightarrow}
\def\vp{\varphi}

\def\vs{\vskip 2em}

\def\v{\vskip 1em}
\def\O{{\cal O}}

\def\C{{\cal C}}

\def\Hat{\widehat}
\def\meas{\hbox{meas}}

\def\dint{\int\!\!\int}

\def\bega{\begin{array}}
\def\enda{\end{array}}
\def\begi{\begin{itemize}}
\def\endi{\end{itemize}}
\def\TV{\hbox{TV}}

\def\ds{\displaystyle}
\def\bel{\begin{equation}\label}
\def\eeq{\end{equation}}
\def\sqr#1#2{\vbox{\hrule height .#2pt
\hbox{\vrule width .#2pt height #1pt \kern #1pt
\vrule width .#2pt}\hrule height .#2pt }}
\def\square{\sqr74}
\def\endproof{\hphantom{MM}\hfill\llap{$\square$}\goodbreak}

\newtheorem{theorem}{Theorem}[section]

\newtheorem{lemma}{Lemma}[section]
\newtheorem{example}{Example}[section]

\newtheorem{definition}{Definition}[section]

\definecolor{bluegreen}{rgb}{0.0, 0.87, 0.87}
\definecolor{blush}{rgb}{0.87, 0.36, 0.51}
\definecolor{bittersweet}{rgb}{1.0, 0.44, 0.37}
\definecolor{burgundy}{rgb}{0.5, 0.0, 0.13}
\definecolor{cardinal}{rgb}{0.77, 0.12, 0.23}

\begin{document}

\title{\bf A Uniqueness Condition for Conservation Laws with Discontinuous Gradient-Dependent Flux}
\vs
\author{Alberto Bressan
and Wen Shen
\\
\, \\
{\small  Department of Mathematics, Penn State University,
University Park, PA 16802, U.S.A. } \,\\ 
\,\\
{\small E-mails:  axb62@psu.edu,~
wxs27@psu.edu}
\,\\   }
\maketitle  

\begin{abstract} 
The paper is concerned with a scalar conservation law with discontinuous gradient-dependent flux.
Namely, the flux is described by two different functions $f(u)$ or $g(u)$, 
when the gradient $u_x$ of the solution is positive or negative, respectively.
In the stable case where $f(u)<g(u)$ for all $u\in\R$, it was proved in \cite{ABS1} that the limits of vanishing viscosity approximations form a contractive semigroup w.r.t.~the $\L^1$ distance.  
Further, they coincide with the limits of  a suitable family of front tracking approximations.
In the present paper we introduce a simple condition that guarantees that every weak, entropy admissible solution of a Cauchy problem 
coincides with the corresponding semigroup trajectory, and hence is unique.
\end{abstract}

\maketitle

\section{Introduction}
\label{sec:1}
\setcounter{equation}{0}

Following the model introduced in \cite{ABS1}, we consider a conservation law  with discontinuous flux, 
where the discontinuity is related to the sign of the gradient $u_x$, namely
\bel{1}
u_t + \Big[\theta(u_x) f(u) + \bigl(1-\theta(u_x) \bigr)g(u)\Big]_x~=~0.\eeq
Here $f,g$ are two different  fluxes, and $\theta$ is a function such that
\bel{2}\theta(s)~=~\left\{\bega{rl} 1\quad &\hbox{if}\quad s>0,\cr
0\quad &\hbox{if}\quad s<0.\enda\right.\eeq
As in \cite{ABS1} we focus on the stable case $f<g$ and assume
\begi
\item[{\bf (A1)}] {\it The flux functions $f,g$ are twice continuously differentiable and satisfy}
\bel{eq:f-le-g}
g(u) - f(u)~\geq ~c_0 ~>~0\qquad \forall u\in\R\,.
\eeq
\endi
To introduce an appropriate definition of weak solution, we observe that for every BV function
$u:\R\mapsto\R$, the distributional derivative $\mu = D_x u$ is a signed Radon measure.
It can be decomposed into a positive and a negative part as $\mu=\mu_+-\mu_-$.
Given a measurable function $\theta:\R\mapsto [0,1]$, the requirement that $\theta=1$ at points where 
$\mu$ is positive, while $\theta=0$ at points where $\mu$ is negative can now be stated
in terms of the measure identity
\bel{mid1}\mu ~=~\chi_{\strut\{ \theta=1\}} \mu_+ - \chi_{\strut\{\theta=0\}} \mu_-\,.\eeq
Here $\chi_{\strut S}$ denotes the characteristic function of a set $S$.
Following \cite{ABS2}, we shall apply the identity (\ref{mid1}) to the measure $D_x u(t,\cdot)$, at a.e.~time $t$.

\begin{definition}\label{d:11} We say that a BV function $u:[0,T]\times \R \mapsto \R$ is a {\bf  weak solution}
to (\ref{1})-(\ref{2}) with initial data
\bel{idata} u(0,x)\,= \,\bar u(x)\eeq
 if the following holds.
\begi
\item[(i)] The map $t\mapsto u(t,\cdot)$ is continuous from $[0,T]$ into $\L^1_{loc}$ and satisfies
the initial condition (\ref{idata}). 
\item[(ii)] There exists a measurable function 
$\theta: [0,T]\times \R\mapsto [0,1]$ such that
\bel{bta} D_xu(t,\cdot)~=~\chi_{\{\theta=1\}}  \Big[D_xu(t,\cdot)\Big]_+ -\chi_{\{\theta=0\}}  \Big[D_x u(t,\cdot)\Big]_-\qquad\hbox{for a.e.}~t\in [0,T],\eeq
and such that, for every compactly supported test function $\vp\in \C^1_c\bigl( ]0,T[\,\times\R\bigr)$,
\bel{wsol} 
\int_0^T \int  \left\{u\vp_t  + \Big[\theta\, f(u) + \bigl(1-\theta\, \bigr)g(u)\Big]\, \vp_x\right\}\, dx\, dt ~=~0.
\eeq
\endi
\end{definition}

Since the integral in (\ref{wsol}) does not depend on the values of $\theta$ on a set of measure zero,
instead of (\ref{bta}) one could simply require that 
$$\theta(t,x)~=~\left\{ \bega{cl} 1\quad \hbox{if}~~u_x(t,x)>0,\\[1mm]
0\quad \hbox{if}~~u_x(t,x)<0.\enda\right.$$
On the other hand, the condition (\ref{bta})  specifies the value of $\theta(t,\cdot)$ 
not only at points where the gradient $u_x$ exists, but also at points where $u(t,\cdot)$ has a jump, 
and on the Cantor part
of the distributional derivative $D_x u(t,\cdot)$. 
This additional information will play a  central  role in the uniqueness
result, which is the main contribution of the present paper. 

Throughout the following, we consider weak solutions having the   additional key property
\begi
\item[{\bf (A2)}] {\it For every $t>0$, the  function $x\mapsto \theta(t,x)$ is continuous.}
\endi

To write the corresponding Rankine-Hugoniot and the Liu admissibility conditions, we recall that a point 
$(\bar t, \bar x)$ is called  a point of {\bf approximate jump} for the function $u=u(t,x)$ if there
exist values $u^-, u^+\in \R$ and a speed $\lambda$ such that the following holds.   Setting
\bel{UT}U(t,x) ~\doteq~\left\{\bega{rl}u^-\quad\hbox{if}~~x<\lambda t,\\
u^+\quad\hbox{if}~~x>\lambda t,\enda\right.\eeq
one has 
\bel{UJ} 
 \lim_{r\to 0+}~\frac{1}{ r^2} \int_{\bar t-r}^{\bar t+r} \int_{\bar x-r}^{\bar x+r} \bigl| u(t,x)- U(t-\bar t, x-\bar x)\bigr|\, dx\, dt~=~0.\eeq
If (\ref{UT})-(\ref{UJ}) hold with $u^-=u^+$ (and hence $\lambda$ arbitrary), we say that $u$ is {\bf approximately 
continuous} at $\bar t, \bar x)$,

In the above setting, if $u,\theta$ provide a weak solution to (\ref{1})-(\ref{2}), a standard argument \cite{Bbook}
yields that, at every point of approximate jump, the left and the right values $u^\pm$ and the speed $\lambda$ in 
(\ref{UT}) must satisfy the Rankine-Hugoniot conditions
\bel{RH} \lambda~=~\frac{\bigl[\theta f(u^+) + (1-\theta) g(u^+) \bigr] - \bigl[\theta f(u^-) + (1-\theta) g(u^-) \bigr]
}{ u^+-u^-}\,,\eeq
where $$\theta~=~\theta(\bar t,\bar x)~=~\left\{ \bega{rl} 1\quad &\hbox{if}~~u^+>u^-,\\[1mm]
 0\quad &\hbox{if}~~u^+<u^-.\enda\right.$$
In addition, we say that the jump satisfies the {\bf Liu admissibility conditions} \cite{Liu} if the following implication hold:
\bel{Liuad} \bega{rl} 
u^-<u^+&\implies\quad 
\displaystyle 
\frac{ f(u^*) - f(u^-) }{ u^*-u^-}~\geq~  \frac{ f(u^+)  - f(u^-) }{ u^+-u^-} \qquad\forall u^*\in \,]u^-, u^+],
\\[4mm]
u^+<u^- &\implies\quad 
\displaystyle
\frac{ g(u^*)  - g(u^-) }{ u^*-u^-}~\geq~  \frac{ g(u^+)  - g(u^-) }{ u^+-u^-} \qquad\forall u^*\in \,]u^+, u^-].\enda
\eeq
More generally, we say that a weak solution $u=u(t,x)$ is {\bf Liu-admissible} if it satisfies the Liu conditions at every point of approximate jump.
In the recent paper \cite{ABS1}, given an initial datum (\ref{idata}), a
solution of  (\ref{1})-(\ref{2}) 
was  constructed as the unique limit of a family of front tracking approximations.
These limit solutions form a semigroup which is contractive w.r.t.~the $\L^1$ distance. 
Furthermore, they also coincide with the unique limits of
parabolic approximations, as $\ve, \delta \to 0$,
\bel{3} u_t +  \Big[\theta_\ve(u_x) f(u) + \bigl(1-\theta_\ve(u_x) \bigr)g(u)\Big]_x
~=~ \delta u_{xx}\,.
\eeq
Here $\theta_\ve(s)~\doteq~\theta(s/\ve)$, while $\theta:\R\mapsto [0,1]$ is a  smooth, nondecreasing function satisfying
\bel{tprop1} \theta(s)~=~\left\{\bega{rl} 1\quad &\hbox{if}\quad s\geq 1,\cr
0\quad &\hbox{if}\quad s\leq -1.\enda\right.
\eeq

 However, the uniqueness of solutions remained an open problem in \cite{ABS1}. The heart of the matter is to
 identify a set of conditions that is satisfied by the semigroup trajectory and by no other solution.  
 Uniqueness results, as in \cite{BDL,BGo, BGu}, become relevant when solutions are constructed as limits  of different types of approximations (relaxation, difference schemes, etc...). In this case, a uniqueness theorem will guarantee that these limit solutions coincide with the unique semigroup trajectory.
 
 As shown by the following counterexample (see Example 1.1 in \cite{ABS1}), 
the Cauchy problem (\ref{1}), (\ref{idata}) can have 
multiple weak solutions, all satisfying the Liu  admissibility conditions.

\begin{figure}[ht]
\centerline{\hbox{\includegraphics[width=15cm]{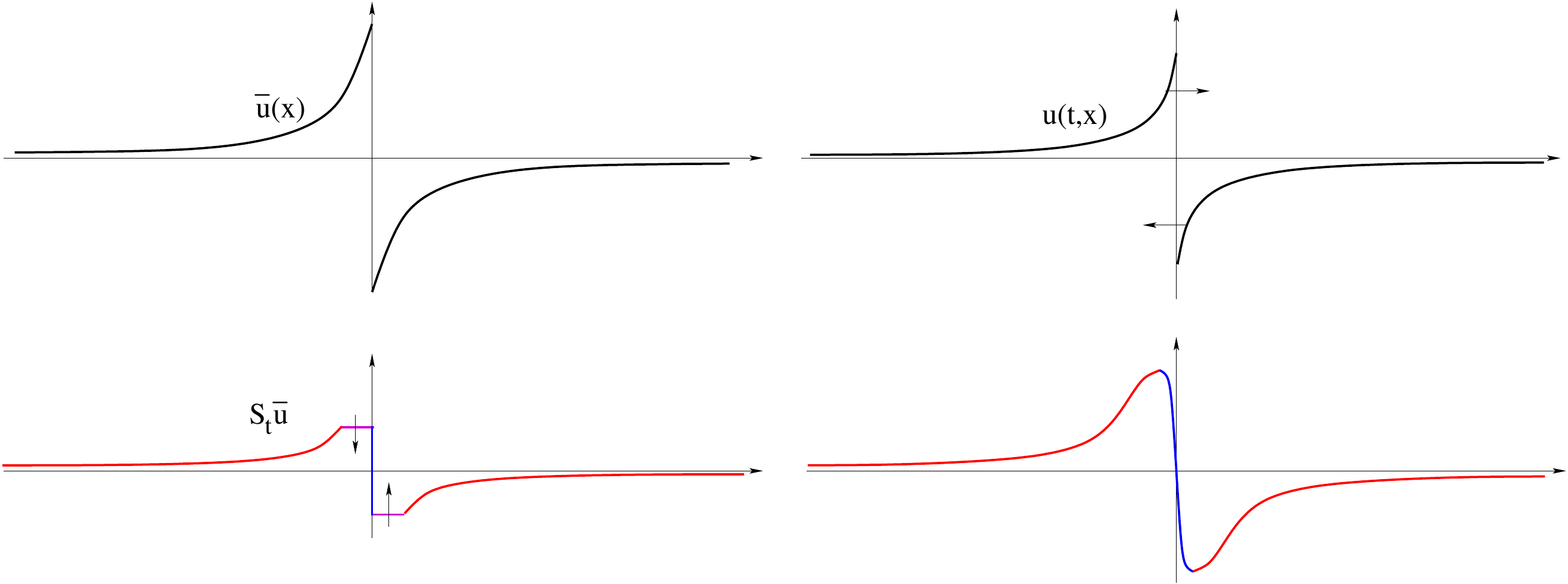}}}
\caption{\small Upper left: the initial data at (\ref{ex2id}). Upper right: the corresponding solution to 
Burgers' equation.   Lower left: the semigroup solution, obtained as limit of front tracking approximations.
Lower right: viscous approximations.  For this example we have $u(t,\cdot)\not= S_t \bar u$ for all $t>0$.
}
\label{f:df106}
\end{figure}

\begin{example}
{\rm 
As flux functions, take
$$f(u)~=~\frac{u^2}{ 2}\,,\qquad\qquad g(u)~=~\frac{u^2}{ 2}+1,$$
and consider the initial data
\bel{ex2id} u(0,x)\,=\,\left\{ \bega{cl} e^x\quad &\hbox{if}\quad x<0,\\
-e^{-x}\quad &\hbox{if}\quad x>0.\enda\right.\eeq
As shown in
Fig.~\ref{f:df106}, let $u=u(t,x)$ be the unique entropy weak solution to Burgers' equation 
$$u_t+f(u)_x~=~0.\qquad\qquad $$
Setting
$$\theta(t,x)~=~\left\{\bega{rl} 1\qquad\hbox{for}~~x\not= 0,\\
0\qquad\hbox{for}~~x= 0,\enda\right.
$$
we obtain
$$
\dint \Big\{u \phi_t + \bigl[ \theta f(u) + (1-\theta) g(u)\bigr] \phi_x\Big\}\, dx\, dt~=~0.
$$

In addition, the downward jump in $u(t,\cdot)$ at $x=0$ clearly satisfies the Liu-entropy admissibility condition, which in this case coincides with  the Lax condition:   $u^->u^+$. }
\end{example}

However,  $u$  in \eqref{ex2id} is  NOT a trajectory of the semigroup.
This can be checked by observing that, for every $t>0$, 
the function  $u(t,\cdot)$ is
strictly increasing on the two open domains where $-\infty<x<0$ and $0<x<\infty$, respectively.
On the other hand, as proved in \cite{ABS1},
the semigroup solution $S_t \bar u$ attains  its maximum value  and its minimum value on two intervals
of positive length.   
These two properties are clearly incompatible.

\begin{figure}[ht]
\centerline{\hbox{\includegraphics[width=15cm]{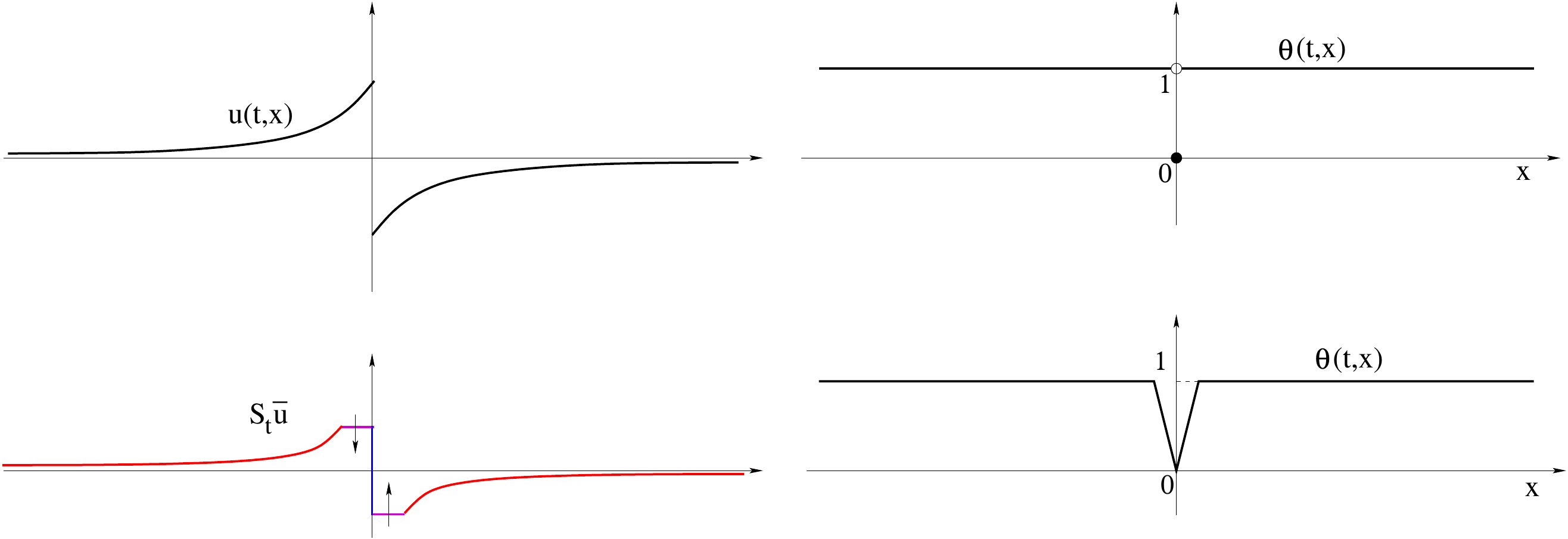}}}
\caption{\small Upper left: the solution to Burgers' equation with initial data (\ref{ex2id}).
Lower left: the semigroup solution. 
Upper and lower right: the corresponding functions $\theta(t,\cdot)$.}
\label{f:df78}
\end{figure}

The corresponding functions $\theta(t,x)$ are shown in Fig.~\ref{f:df78} (right),
for the solution of Burgers' equation and for the semigroup solution.   
A major difference between them is the following. 
In the solution of Burgers' equation, $\theta(t,\cdot)$ has a downward jump at the origin, 
for every $t\geq 0$.
In the semigroup solution the function $\theta$ is continuous for all $t>0$.

%
Purpose of the present note is to prove that the continuity assumption {\bf (A2)} precisely characterizes the 
vanishing viscosity limits.   Namely, any entropy admissible weak solution $u$ of the Cauchy problem (\ref{1}), (\ref{idata}),  with $\theta=\theta(t,x)$ continuous
in the $x$ variable, is unique and coincides with the corresponding semigroup trajectory.  

 The remainder of the paper is organized as follows. 
A precise statement of the result is given in 
Section~\ref{sec:2}.  Details of the proof are then worked out in Section~\ref{sec:3}.

For a  general introduction to conservation laws and admissibility conditions, we refer to 
\cite{Bbook, Cbook, D2, HR, Lbook}.

\section{A uniqueness theorem}
\label{sec:2}
\setcounter{equation}{0}
To simplify some of the analysis,
as in \cite{ABS1}, we consider the space $\L^1_{per}(\R)$ of locally integrable, periodic functions, so that $w(x)=w(x+1)$ for all $x\in \R$.    The $\L^1$ norm and the total variation are thus 
computed over the domain 
$[0,1]$
with endpoints identified.
By the construction of solutions to (\ref{1}) as limits of front tracking approximations \cite{ABS1}, it is established that semigroup solutions $u(t,\cdot)=S_t\bar u$ have the following properties.
\begi
\item[(i)] The map $t\mapsto u(t,\cdot)=S_t{\bar u} \in \L^1_{per}$ is uniformly Lipschitz continuous on every  
interval $[t_0,T]$ with $t_0>0$.
\item[(ii)] Every semigroup trajectory is a weak solution according to Definition~\ref{d:11}, and satisfies 
the Liu admissibility 
conditions (\ref{Liuad}) at every point of approximate jump.
\item[(iii)] For each fixed $t>0$, the solution $u(t,\cdot)$ has bounded variation. 
It attains its maxima and minima on a finite
number of nontrivial disjoint intervals (over one spatial period).  
The number of these distinct intervals is non-increasing in time.
The total variation $\TV\bigl\{ u(t,\cdot)\bigr\}$ 
is strictly decreasing in time, until the solution becomes constant.
\item[(iv)] For each $t>0$ the corresponding function $\theta(t,\cdot)$ is piecewise affine. It increases from 0 to 1
on intervals where $u(t,\cdot)$ attains a local minimum, and decreases from 1 to 0 on intervals where 
$u(t,\cdot)$ attains a local maximum. 
\endi

\begin{figure}[ht]
\centerline{\hbox{\includegraphics[width=10cm]{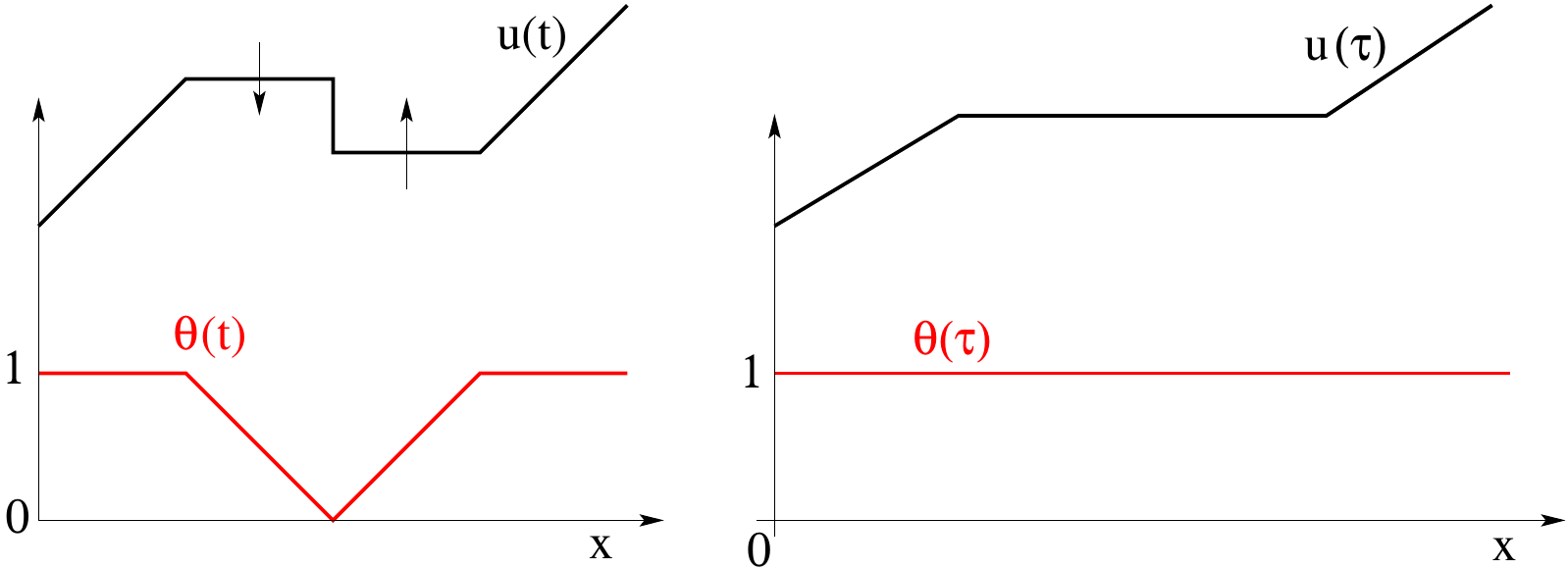}}}
\caption{\small An example where, for the semigroup solution $u(t)= S_t \bar u$, the corresponding function $\theta(t,x)$ is not continuous w.r.t.~time.  For every time  $t>0$ the function $\theta(t,\cdot)$ is piecewise affine.
If $\tau$ is the time when the two horizontal portions of the graph of $u(t,\cdot)$ join, 
the function $\theta(\tau,\cdot)$ does not coincide with the limit of $\theta(t,\cdot)$ as $t\to \tau-$.}
\label{f:df79}
\end{figure}

It is worth noting that, even for semigroup trajectories, the map $\theta(t,x)$ can be discontinuous w.r.t.~time.
An example where this happens is when $u(t,\cdot) $ attains a local maximum and a local minimum on two adjacent
intervals. In finite time these intervals
merge together,  as shown in Fig.~\ref{f:df79}.

The following uniqueness theorem provides a characterization of  the semigroup  solutions of (\ref{1})-(\ref{2}).

\begin{theorem} \label{THM}
Let $f,g$ satisfy the assumptions {\bf (A1)}. Given an initial data $\bar u\in \L^1_{per}$, let $u=u(t,x)$ be a 
Liu admissible weak solution of 
the Cauchy problem (\ref{1})-(\ref{2}) with initial data  (\ref{idata}). Assume that
\begi
\item[(i)] Given  any $t_0>0$, for $t\in [t_0,T]$ the function $u(t,\cdot)$ 
and the corresponding function $\theta(t,\cdot)$ have uniformly
bounded variation.
\item[(ii)] For every $t>0$ the  function $\theta(t,\cdot)$ is continuous.
\endi
Then one has \bel{usem} u(t,\cdot)\,=\, S_t{\bar u}\qquad\forall t\in [0,T].\eeq In other words, the solution 
with the above properties is unique and coincides with the corresponding semigroup trajectory.
\end{theorem}

\section{Proof of Theorem \ref{THM}}
\label{sec:3}
\setcounter{equation}{0}
 In this section we prove Theorem~\ref{THM}  in several steps.

{\bf 1.} For any $t_0>0$, 
the assumptions
\bel{utBV}
\TV \bigl\{u(t,\cdot)\bigr\}\,\leq\,C_0\qquad   \TV \bigl\{\theta (t,\cdot)\bigr\}\,\leq\,C_0
\qquad \forall t\in [t_0,T]
\eeq
yield a bound on the total variation of the flux function
\bel{Flux}x~\mapsto~ F(t,x)~\doteq~\theta(t,x) f\bigl(u(t,x)\bigr) + \bigl(1- \theta(t,x)\bigr) g\bigl(u(t,x)\bigr) \eeq
in terms of the constant $C_0$ and of the Lipschitz constants of $f$ and $g$.
As a consequence (see Theorem~4.3.1 in \cite {D2} for details),
restricted to 
the time interval $[t_0, T]$,
 the map
$t\mapsto u(t,\cdot)$ is Lipschitz continuous with values in $\L^1_{per}$. 
For some constant $L_0$ we thus have
\bel{L0} \bigl\| u(t_1,\cdot)- u(t_2,\cdot)\bigr\|_{\L^1_{per}}~\leq~L_0 \, |t_1-t_2|\qquad\qquad \forall t_1, t_2\in [t_0, T].
\eeq

Since $u=u(t,x)$ is a BV function of two variables on the domain $[t_0,T]\times [0,1]$,
by a well known structure theorem \cite{AFP, EG}
there exists a set of times ${\cal N}$ of measure zero such that the following holds.
For every $t\in [t_0,T]\setminus {\cal N}$ and every $x\in [0,1]$, the point $(t,x)$ is either a point of approximate continuity, or a point of approximate jump for $u$.

\v
{\bf 2.} The proof will be achieved by showing that, for any $t_0>0$, 
\bel{use2}
u(t)~=~S_{t-t_0} u(t_0),\qquad\qquad \forall t\in [t_0,T]\,.\eeq
By the contractivity of the semigroup, this will imply
\bel{tine}
\bigl\|u(t,\cdot)-S_t \bar u\bigr\|_{\L^1_{per}}~=~\bigl\|S_{t-t_0}u(t_0)-S_t \bar u\bigr\|_{\L^1_{per}}\\[2mm]
~\leq~\bigl\|u(t_0)-S_{t_0} \bar u\bigr\|_{\L^1_{per}}\,.\eeq
By the continuity of the maps $t\mapsto u(t,\cdot)$ and $t\mapsto S_t\bar u$, the right hand side goes to zero as 
$t_0\to 0$. This will yield (\ref{usem}), proving uniqueness.

In the remaining steps we will prove the identity (\ref{use2}).  Thanks to the Lipschitz continuity of the map $t\mapsto u(t,\cdot)\in \L^1_{per}$ and the contractivity of the semigroup, 
the following
elementary error estimate holds (see \cite{Bbook})
\bel{ees} \bigl\|u(t)-S_{t-t_0} u(t_0)\bigr\|_{\L^1_{per}} ~\leq~\int_{t_0}^t
\liminf_{h\to 0+} ~\frac{1}{ h} \Big\| u(\tau+h) - S_h u(\tau)\Big\|_{\L^1_{per}}\, d\tau\,.\eeq

\begin{figure}[ht]
\centerline{\hbox{\includegraphics[width=7.5cm]{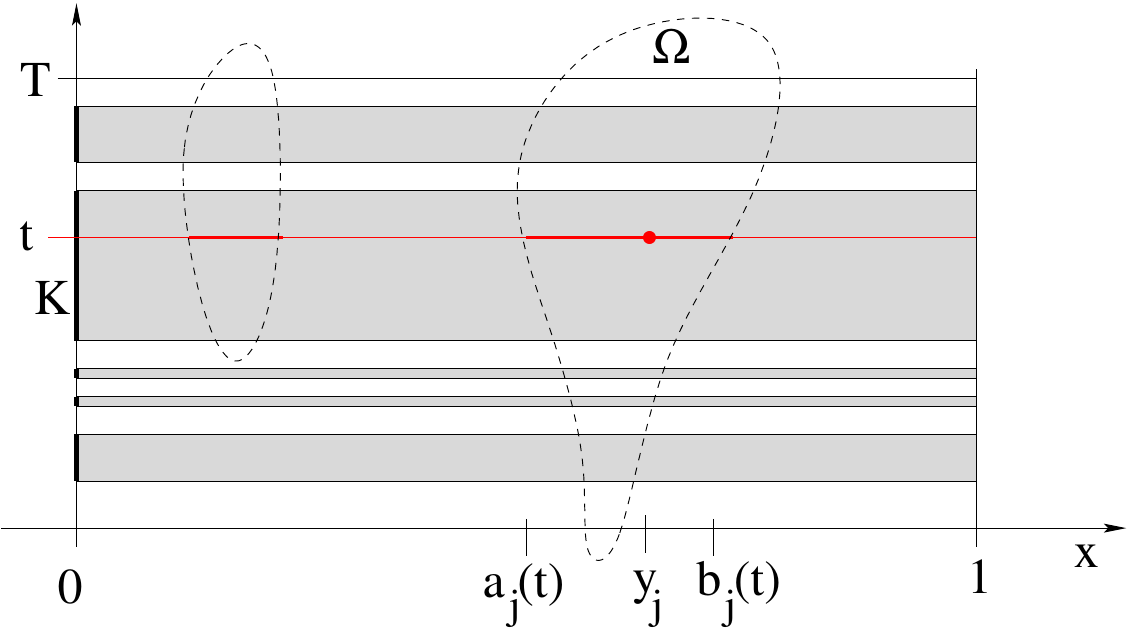}}}
\caption{\small The compact set $K\subseteq [0,T]$  
and $\Omega\subset \R^2$, considered at (\ref{Om}).}
\label{f:df85}
\end{figure}

\v
{\bf 3.} By assumption, the function $\theta=\theta(t,x)$ is measurable in $t$ and continuous in $x$.
By a theorem of Scorza-Dragoni \cite{SD}, for any $\ve_0>0$ there exists a compact set $K\subseteq [0,T]$
with $\meas(K)>T-\ve_0$ such that $\theta$ is continuous (w.r.t.~both variables $t,x$) restricted to the product space $K\times [0,1]$.
Within the set $K\times [0,1]$ we distinguish
\begi
\item[(i)] a compact region where $\theta=0$,
\item[(ii)]  a compact region where $\theta=1$,
\item[(iii)] a relatively open region where $0<\theta<1$.
\endi
As shown in Fig.~\ref{f:df85},
we can find an open set $\Omega\subset\R^2$ such that 
\bel{Om}
\Big\{ (t,x)\in K\times [0,1]\,;~~0<\theta(t,x)<1\Big\}~=~\bigl(K\times [0,1]\bigr)\cap\Omega.\eeq

By the continuity of $\theta$, for every $t\in K$ the set where $0<\theta(t,\cdot)<1$ is a countable union of  open intervals
$\,\bigl]a_i(t), b_i(t)\bigr[\,$, $i\geq 1$. 
\v
{\bf 4.} For any rational point $y_j\in [0,1]$, call 
$\,\bigl]a_j(t), b_j(t)\bigr[\,$ the connected component of the set
$$\bigl\{ x\in \,]0,1[\,;~~(t,x)\in \Omega\bigr\}$$
which contains $y_j$.  Since $\Omega\subset \R^2$ is open,  this interval is well defined for all $t$ 
such that $(t, y_j)\in \Omega$.

The condition (\ref{bta}) implies that the function 
$u(t,\cdot)$ is constant on each of these open intervals, say
\bel{uuj} u(t,x)~=~u_j(t)\qquad\forall x\in\,\bigl]a_j(t), b_j(t)\bigr[\,,\quad t\in K.\eeq
We claim that the map $t\mapsto u_j(t)$ is locally Lipschitz continuous, restricted to $t\in K$.

Indeed, 
 for any $(\tau,y_j)\in\Omega$,  there exists $\delta,\rho>0$ such that 
$$b_j(t)-a_j(t)\geq\rho\qquad\qquad \hbox{whenever} ~~~|t-\tau|\leq \delta.$$
Recalling (\ref{L0}) 
we obtain
$$\rho\cdot \bigl| u_j(t_1)- u_j(t_2)\bigr|~\leq~\bigl\| u(t_1,\cdot) - u(t_2,\cdot )\bigr\|_{\L^1_{per}}~\leq~L_0\cdot |t_1-t_2|\,.$$
This proves Lipschitz continuity, with constant $L_0/\rho$.  

Observing that the complement $[0,T]\setminus K$ is a countable union of open intervals $\,]t_{k-1}, t_k[\,$, $k\geq 1$,
we can extend the function $u_j(\cdot)$ in an affine way on each of these intervals.
This yields a function $t\mapsto u_j(t)$ defined on the entire interval $[\tau-\delta, \tau+\delta]$,
which is still   Lipschitz continuous with the same Lipschitz constant 
$L_0/\rho$.

 As a consequence, by Rademacher's theorem
the time derivative $\dot u_j(t)$ exists at a.e.~time $t$ such that  $(t,y_j)\in\Omega$.
In particular, the set of times
\bel{Ndef}
{\cal N}'~\doteq~\Big\{ t\in K\,;
~~(t, y_j) \in \Omega ~~\hbox{for some}~j\geq 1, ~~~u_j (\cdot)~\hbox{is not differentiable at}~t\Big\}\eeq
has  Lebesgue measure zero.

%

%
\v
{\bf 5.} Let $\tau$ be a Lebesgue point of the set $K$, so that
\bel{lebpo}  \tau\in K,\qquad 
\lim_{h\to 0+} \frac{1}{ h} \meas\Big( [\tau-h, \tau+h]\setminus K\Big)~=~0.\eeq
Assume $\tau\notin {\cal N}'$.
If $y_j$ is any rational point such that $(\tau,y_j)\in\Omega$, then the open intervals 
$\,\bigl]a_j(t), b_j(t)\bigr[\,$ are well defined for all $t$ in a neighborhood of $\tau$.
Let the function $u_j(t)$ be as in (\ref{uuj}).   By assumption, the derivative $\dot u_j(\tau)$ exists.

We claim that,   on the interval $\bigl]a_j(\tau),b_j(\tau)\bigr[ \,$,
one of the following alternatives holds:
\bel{ujt}\bega{ll}
\dot u_j(\tau) =\displaystyle  
\frac{ f(u_j(\tau))- g(u_j(\tau))}{ b_j(\tau)-a_j(\tau)}\,<\,0\quad 
&\hbox{if} \quad \theta(\tau, a_j(\tau))=1, ~\theta(\tau, b_j(\tau))=0,\\[4mm]
\dot u_j(\tau)=\displaystyle 
\frac{ g(u_j(\tau))- f(u_j(\tau))}{ b_j(\tau)-a_j(\tau)}\,>\,0\quad &\hbox{if} \quad \theta(\tau, a_j(\tau))=0, ~\theta(\tau, b_j(\tau))=1.\enda
\eeq
Indeed, fix $x_1, x_2 \in \bigl]a_j(\tau), b_j(\tau)\bigr[\,$ with $x_1<x_2$,
and let $h>0$ be small enough.  
Applying the conservation law over the rectangle $[x_1,x_2]\times [\tau, \tau+h]$, 
by the divergence theorem we obtain
\begin{align*}
&\hspace{-5mm} \bigl[u_j(\tau+h)-u_j(\tau)\bigr] \cdot (x_2-x_1) 
\\
&=~ \int_{\tau}^{\tau +h}  \Big[ \theta(t, x_1) f(u_j(t)) + (1-\theta(t,x_1)) g(u_j(t)) 
\Big] dt 
\\
&\quad  -  \int_{\tau}^{\tau +h}  \Big[ \theta(t,x_2) f(u_j(t)) + (1-\theta(t,x_2)) g(u_j(t)) 
\Big] dt 
\\
&=~ \int_{\tau}^{\tau +h}  \bigl[ \theta(t, x_1) -  \theta(t, x_2) \bigr] \cdot 
\bigl[ f(u_j(t)) - g(u_j(t)) \bigr] dt .
\end{align*}
%
Taking the limit as $h\to 0$, using the fact that $\theta$ is continuous restricted to $K\times [x_1, x_2]$ and 
$\tau$ is a Lebesgue point of $K$,  we obtain
$$
\dot u_j(\tau) =\frac{ \theta(\tau, x_1) -  \theta(\tau, x_2) }{x_2-x_1} \cdot 
\bigl[ f(u_j(\tau)) - g(u_j(\tau)) \bigr] .
$$
Since the left hand side does not depend on the space variable $x$, 
the quotient
$$ 
\frac{ \theta(\tau, x_1) -  \theta(\tau, x_2) }{x_2-x_1}
$$ 
must be a constant,  independent of the  choice of $x_1,x_2$. 
We conclude that  $\theta$ must be an affine function, taking values in the set $\{0,1\}$ at the endpoints
of the interval $\bigl[a_j(\tau), b_j(\tau)\bigr]$.   By assumption, $\theta(\tau,x)\notin \{0,1\}$
for $a_j(\tau)< x <b_j(\tau)$, hence $\theta(\tau,\cdot)$ cannot be constant.
Therefore \eqref{ujt} holds. 
Letting $x_1 \to a_j(\tau) $ and $ x_2 \to b_j(\tau)$, we thus obtain
$$\theta_x(\tau,x)~=~\hbox{constant}~=~\frac{\pm 1}{ b_j(\tau)-a_j(\tau)}
\qquad\qquad \hbox{for}~~x\in \bigl]a_j(\tau), b_j(\tau)\bigr[\,.$$
This proves our claim.
\v
{\bf 6.} Next, consider any time
$\tau\in K$.
By the previous step, on each interval $\,\bigl]a_j(\tau), b_j(\tau)\bigr[\,\subset [0,1]$ where $0<\theta(\tau,x)<1$, the 
function $\theta(\tau,\cdot)$ is affine, taking values 0 and 1 at the two endpoints.  
By the continuity  of $\theta(\tau,\cdot)$, there can be only finitely many 
such intervals.  After a relabeling, we can assume these intervals are $\,\bigl]a_i(\tau), b_i(\tau)\bigr[\,$,
$i=1,\ldots,N(\tau)$.

Since $\theta$ is uniformly continuous restricted to $K\times [0,1]$, 
and  $\theta(t,\cdot)$
affine on each $\,\bigl]a_i(t), b_i(t)\bigr[\,$ where it varies between  0 and 1,
for some constant $M$ we must have
\bel{txb}
\bigl|\theta_x(t,x)\bigr|~\leq~M\qquad\qquad \forall t\in K\quad \hbox{and for a.e.~} x\in [0,1].\eeq
In turn, this implies that all the maps
$$t\mapsto a_i(t),\qquad\qquad t\mapsto b_i(t)$$
are uniformly continuous restricted to $K$.
By choosing $\bar h>0$ sufficiently small, we conclude that, for each $t\in [\tau, \tau+\bar h]\cap K$,
the function $\theta(t,\cdot)$ is piecewise affine, and varies between 0 and 1 on a fixed number of
non-overlapping
intervals $\bigl[a_i(t), b_i(t)\bigr]$, $i=1,\ldots,N$.  Moreover, the endpoints $a_i, b_i$ are continuous functions 
of time, restricted to $K$.

 On each interval $\bigl[b_i(t), a_{i+1}(t)\bigr]$,  the function $\theta(\tau,\cdot)\in\{0,1\}$ is constant. 
We further distinguish 5 cases.
\begi
\item[(1)] $b_{i}(t)= a_{i+1}(t)$;
\item[(2)]  $b_{i}(t)< a_{i+1}(t)$ and $u(t,\cdot)$ has a jump at both $b_{i}(t)$ and $a_{i+1}(t)$;
\item[(3)]  $b_{i}(t)< a_{i+1}(t)$ and $u(t,\cdot)$ is continuous at both $b_{i}(t)$ and $a_{i+1}(t)$;
\item[(4)] $b_{i}(t)< a_{i+1}(t)$ and $u(t,\cdot)$ is continuous at $b_{i}(t)$ but has a jump at $a_{i+1}(t)$;
\item[(5)] $b_{i}(t)< a_{i+1}(t)$ and $u(t,\cdot)$ has a jump at $b_{i}(t)$ but is continuous at $a_{i+1}(t)$.
\endi

\begin{figure}[ht]
\centerline{\hbox{\includegraphics[width=9cm]{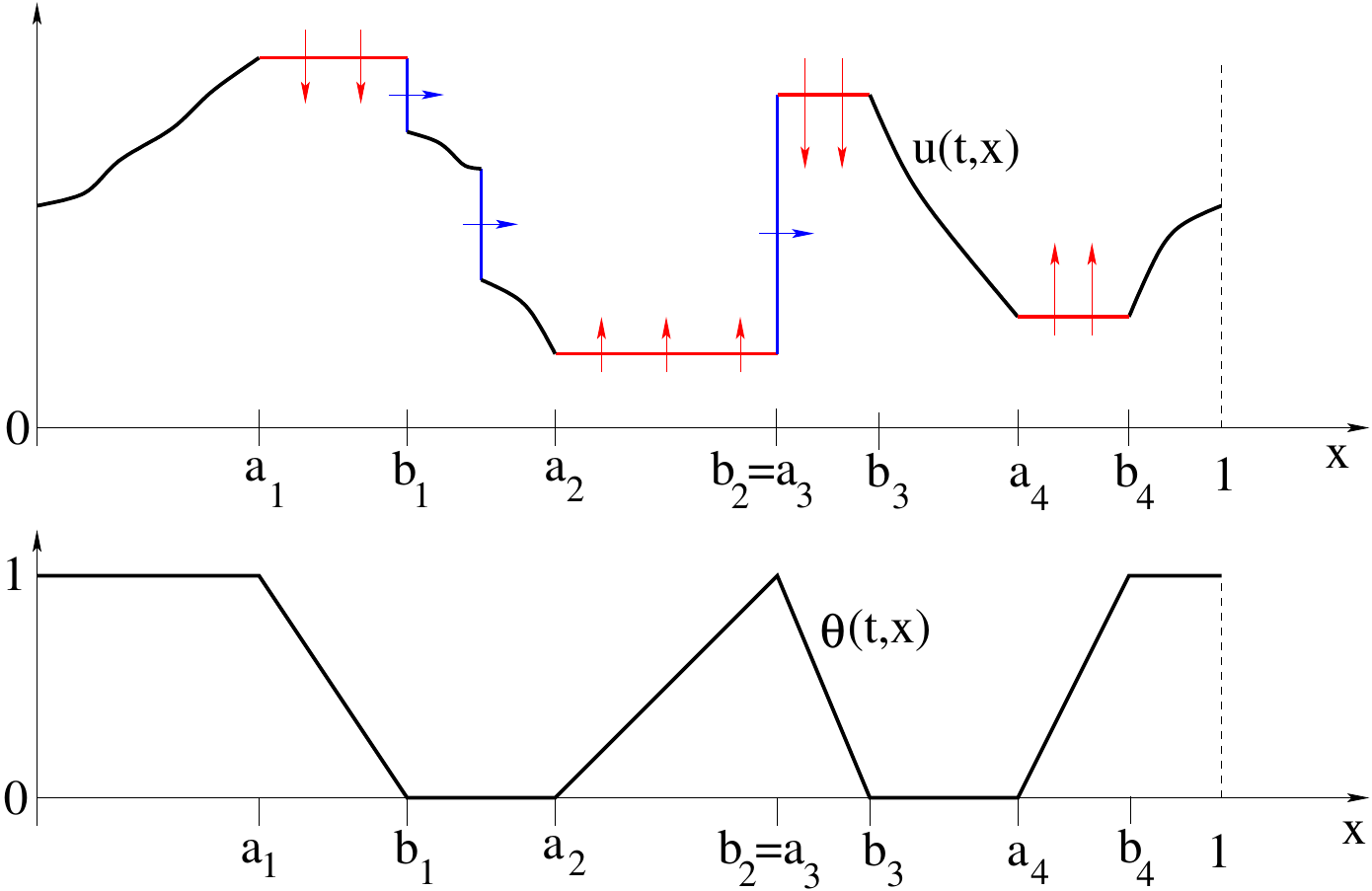}}}
\caption{\small 
Top: a BV solution to the conservation law (\ref{1})-(\ref{2}), at a given time $t>0$.   Bottom: the corresponding function 
$\theta$. Note that $u(t,\cdot)$ is constant on every  open interval $\,\bigl]a_i(t), b_i(t)\bigr[\,$. 
}
\label{f:df95}
\end{figure}

Some of these cases are illustrated in Figure~\ref{f:df95}.

We can then split the set $K$ into $5^N$ measurable subsets $K_\alpha(t)$, distinguishing 
the 5 possible cases for every interval $[b_i(t), a_{i+1}(t)\bigr]$, $i=1,\ldots,N$.
More precisely, 
given a multi-index
$\alpha=(\alpha_1,\ldots,\alpha_n)\in \{1,2,3,4,5\}^N$, we thus define the set of times
\bel{Ka}
K_\alpha \,\doteq\, \Big\{ t\in K\,;   ~\hbox{Case $\alpha_i$ occurs for the interval}~ 
\bigl[b_i(t), a_{i+1}(t)\bigr],~~i=1,\ldots, N 
\Big\}.
\eeq
Clearly the above sets are non-overlapping,
$\bigcup_\alpha K_\alpha=K$ and a.e.~time $t\in K$ is a Lebesgue point of one of the sets $K_\alpha$.
\v

{\bf 7.}  After these preliminaries, the heart of the proof consists in showing that
\bel{key}
\liminf_{t\to \tau+} ~\frac{1}{t-\tau} \Big\| u(t) - S_{t-\tau} u(\tau)\Big\|_{\L^1_{per}}
~=~0
\qquad\quad\hbox{for a.e.}~~\tau\in K.
\eeq
If (\ref{key}) holds, recalling that  both the solution $t\mapsto u(t)$ and the semigroup trajectory $t\mapsto S_{t-t_0} u(t_0)$ are Lipschitz continuous with values in $\L^1_{per}$ for some Lipschitz constant $L_0$, from (\ref{ees}) we immediately
obtain 
\bel{ees2} 
\bega{l}\ds \bigl\|u(t)-S_{t-t_0} u(t_0)\bigr\|_{\L^1_{per}} 
~\leq~
\int_{t_0}^t 
\liminf_{s\to \tau+} ~\frac{1}{s-\tau+} \Big\| u(s) - S_{s-\tau} u(\tau)\Big\|_{\L^1_{per}}
\, d\tau\\[4mm]
\qquad \ds\leq~\int_{[t_0, t]\setminus K} 
2 L_0\, d\tau~=~2 L_0\cdot \meas\Big([t_0, t]\setminus K\Big)~\leq~ 2L_0\ve_0
\,.\enda \eeq
Since $\ve_0>0$ can be chosen arbitrarily small, this will imply (\ref{use2}), proving the theorem.

In the remaining steps of the proof we will establish (\ref{key}). 
The general strategy for proving (\ref{key}) goes as follows.  
Given  $\ve>0$,   for all $\tau\in K$ (outside a set of times of measure zero),
we construct a leading order approximation $U^{app}(t,\cdot)$ with the property that
\bel{key2}\liminf_{t\to\tau+} ~ \frac{1}{ t-\tau} \Big\| u(t) -U^{app}(t)\Big\|_{\L^1_{per}}
~\leq~C\ve,\eeq
for some constant $C$ independent of $\ve$.  Moreover, we show that the same inequality holds
if $u$ is replaced  with any other weak solution $U$ which satisfies the assumption of the theorem and 
coincides with $u$ at time $t=\tau$.   In particular, this holds when $U(t) = S_{t-\tau} u(\tau)$ is the corresponding semigroup trajectory.  In turn,  by the triangle inequality this implies
\bel{key3}\liminf_{t\to\tau+} ~ \frac{1}{ t-\tau} \Big\| u(t) -S_{t-\tau} u(\tau)\Big\|_{\L^1_{per}}
~\leq~2C\ve\eeq
for any $\ve>0$, and hence (\ref{key}).

To establish \eqref{key2}, 
given $\ve>0$ the approximation $U^{app}$ will be constructed by 
covering the domain 
$[0,1]$ with finitely many open intervals $]p_k, q_k[$, $k=1,\ldots,N$. 
We then construct functions $U_k^{app}$ 
on each interval, so that, for some $\lambda^*>0$ large enough,
\bel{key4}\liminf_{t\to\tau+} ~\frac{1}{ t-\tau} \int_{p_k+\lambda^*(t-\tau)}^{q_k-\lambda^*(t-\tau)}
\Big| u(t,x) -U_k^{app}(t,x)\Big|\, dx~\leq~C'\ve,
\qquad k=1,\ldots,N,\eeq
where the constants $N, C'$ are independent of $\ve$.
Since the open intervals $]p_k, q_k[$ cover the entire domain $[0,1]$, the same holds for the intervals of integration
in (\ref{key4}), for all $t$ sufficiently close to $\tau$.
\v
{\bf 8.} To fix ideas, let $\tau$  be a Lebesgue point of a 
set $K_\alpha\subseteq K$, with  $\tau \notin ({\cal N}\cup {\cal N}')$.
Let $\bigl[a_i(\tau), b_i(\tau)\bigr]$, $i=1,\ldots,N$, be the corresponding intervals where 
$\theta(\tau,\cdot)$ varies between 0 and 1, and $u(\tau,\cdot)$ is constant.

\begin{figure}[ht]
\centerline{\hbox{\includegraphics[width=13cm]{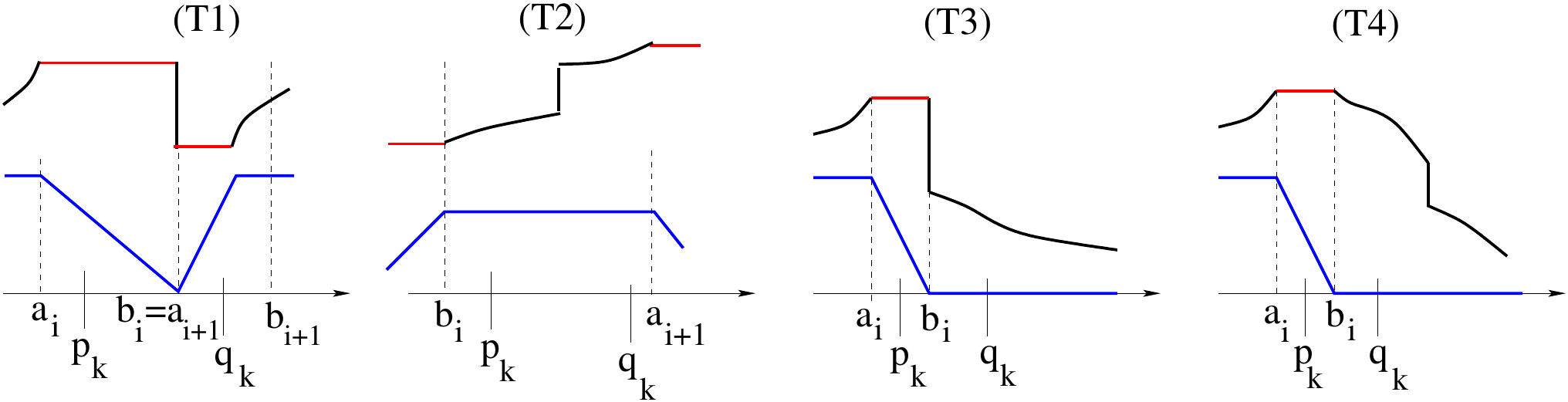}}}
\caption{\small The function $u(\tau,x)$ and the corresponding function $\theta(\tau,x)$, in the configurations (T1)--(T4). }
\label{f:df96}
\end{figure}

Given $\ve>0$, we cover the interval $[0,1]$ with open intervals $\, ]p_k, q_k[\,$ 
such that:
\begin{itemize}
\item Each interval $\, ]p_k, q_k[\,$ contains at most one of the points $a_i(\tau), b_i(\tau)$, $i=1,\ldots,N$.  If $a_{i+1}(\tau)=b_i(\tau)$, it is treated as one single point. 
\item 
If $p_k<b_i(\tau)<q_k$, then  $u(\tau,\cdot)$ is constant on $\bigl]p_k,  b_i(\tau)\bigr[$, 
and 
\bel{tv1}
\TV\Big\{ u(\tau,\cdot)\,;~\bigl]b_i(\tau), q_k\bigr[\,\Big\}\, <\,\ve.
\eeq
Note that \eqref{tv1}  can be achieved by taking  $q_k$ sufficiently close to $b_i(\tau)$.
\item  
If $p_k<a_i(\tau)<q_k$,  then $u(\tau,\cdot)$ is constant on $\bigl]a_i(\tau),  q_k\bigr[$,
and 
\bel{tv1a}
\TV\Big\{ u(\tau,\cdot)\,;~\bigl]p_k,  a_i(\tau)\bigr[\,\Big\}\, <\,\ve. 
\eeq
Note that \eqref{tv1a}  can be achieved by taking  $p_k$ sufficiently close to
 $a_i(\tau)$.

\end{itemize}
We note that at most $3N$ open intervals $\, ]p_k, q_k[\,$ are needed. 
Among them, four different types of intervals shall be considered (see Fig.~\ref{f:df96}).
\begi
\item[(T1)] $p_k<b_i(\tau) = a_{i+1}(\tau)<q_k$;
\item[(T2)] $ b_i(\tau)< p_k< q_k < a_{i+1} (\tau)$;
\item[(T3)] $p_k<b_i(\tau) < q_k< a_{i+1}(\tau)$ and $u(\tau, \cdot)$ has a jump at $x=b_i(\tau)$;
\item[(T4)] $p_k<b_i(\tau) < q_k< a_{i+1}(\tau)$ and $u(\tau, \cdot)$ is continuous at $x=b_i(\tau)$.
\endi
The behavior of the solution on intervals of the types (T1)--(T4) will be separately analyzed in the forthcoming steps.

%

\v

{\bf 9.} We begin by studying the behavior of the solution $u$ on an interval of type (T1).  To fix ideas,
assume that $u(\tau,\cdot)$ is constant for $x\in \bigl[p_k, b_i(\tau)\bigr[\,$ and  $x\in \,\bigl]b_i(\tau), q_k\bigr]$,
and $u(\tau,\cdot)$ has a downward (Liu admissible)
jump at $x=b_j(\tau) = a_{i+1}(\tau)$.

Assuming that $\tau$ is a Lebesgue point of the set $K_\alpha$ of times where this behavior occurs,
we can find a decreasing sequence of times $t_n\in K_\alpha$, $t_n\to \tau$, such that
$b_i(t_n)= a_{i+1}(t_n)$ for every $n\geq 1$.

By the analysis in step {\bf 4}, there exists $\delta>0$ small enough such that, for $t\in [\tau, \tau+\delta]\cap K$,
we have
\begi
\item
$u(t,x)=u_i(t)$  for $p_k\leq  x<b_i(t)$,
\item  $u(t,x)=u_{i+1}(t)$  for $a_{i+1}(t)<  x\leq q_k$.
\endi
Here $b_i, a_{i+1}: K\mapsto\R$ are continuous functions with  $b_i(t) =a_{i+1}(t)$. 
 Meanwhile $u_i, u_{i+1}:K\mapsto \R$ are Lipschitz continuous
functions of time, where $u_i$ is decreasing and  $u_{i+1}$ is increasing.
Moreover, at time $t=\tau$ the time derivatives 
\bel{dotsig} 
\dot u_i(\tau)\,=\, \frac{f(u_i(\tau)) -g(u_i(\tau)) }{ b_i(\tau)-a_i(\tau)} \,<\,0,\quad
 \dot u_{i+1}(\tau)\,= \, \frac{g(u_{i+1}(\tau)) -f(u_{i+1}(\tau)) }{ b_{i+1}(\tau)-a_{i+1}(\tau)} \,>\,0\eeq
are well defined.

Since $\tau\notin {\cal N}$, the solution $u$ has an approximate jump at the point $\bigl(\tau, b_i(\tau)\bigr)$,
with Rankine-Hugoniot speed
\bel{dotb}\lambda_i~\doteq~\frac{g\bigl(u_i(\tau)\bigr) - g\bigl(u_{i+1}(\tau)\bigr) 
}{ u_i(\tau) -u_{i+1}(\tau)}\,.
\eeq
Introducing the function
\bel{Uidef}
U_i(t,x)~\doteq~\left\{\bega{rl} u_i(\tau) + (t-\tau) \dot u_i(\tau) \qquad &\hbox{if}\quad x< b_i(\tau) + (t-\tau)
\lambda_i\,,\\[2mm]
 u_{i+1}(\tau) + (t-\tau)\dot u_{i+1}(\tau) \qquad &\hbox{if}\quad x> b_i(\tau) + (t-\tau) \lambda_i\,,
 \enda\right.\eeq
 by (\ref{dotsig})-(\ref{dotb})
it follows 
\begin{eqnarray}
&& \hspace{-1cm} 
\lim_{t\to \tau+,\, t\in K_\alpha} ~\frac{1}{ t-\tau}\int_{p_k}^{q_k} \Big|u(t,x) - U_i(t,x)\Big|\, dx 
\nonumber\\
&=& \lim_{t\to \tau+,\, t\in K_\alpha} ~\frac{1}{ t-\tau}\int_{p_k}^{b_i(\tau)+(t-\tau)\lambda_i} \Big|u_i(t) - u_i(\tau) -(t-\tau)\, \dot u_i(\tau)\Big|\, dx 
\nonumber\\
&& +\lim_{t\to \tau+,\, t\in K_\alpha} ~\frac{1}{ t-\tau}\int_{b_i(\tau)+(t-\tau)\lambda_i}^{q_k} \Big|u_{i+1}(t) - u_{i+1}(\tau) -(t-\tau)\, \dot u_{i+1}(\tau) \Big|\, dx
\nonumber\\
&=& \lim_{t\to \tau+,\, t\in K_\alpha} ~[b_i(\tau)+(t-\tau)\lambda_i-p_k ] \cdot
\left| \frac{u_i(t) - u_i(\tau)}{t-\tau} - \dot u_i(\tau)\right|\, 
\nonumber\\
&& + \lim_{t\to \tau+,\, t\in K_\alpha} ~[q_k - b_i(\tau)-(t-\tau)\lambda_i ] \cdot
\left| \frac{u_{i+1}(t) - u_{i+1}(\tau)}{t-\tau} - \dot u_{i+1}(\tau)\right| 
\nonumber\\
&=& 0\,.
\label{lim4} 
\end{eqnarray}
On the other hand, denoting by $U(t,\cdot)= S_{t-\tau} u(\tau)$ the semigroup trajectory, we also have
\bel{lim44} 
 \lim_{ t\to \tau+}
  ~\frac{1}{ t-\tau}\int_{p_k}^{q_k} \Big|U(t,x) - U_i(t,x)\Big|\, dx~=~0.\eeq
Together, (\ref{lim4})-(\ref{lim44}) imply
\bel{lim45}  \lim_{ t\to \tau+,\, t\in K_\alpha}  ~\frac{1}{ t-\tau}\int_{p_k}^{q_k} \Big|u(t,x) - U(t,x)\Big|\, dx~=~0.\eeq

\v
{\bf 10.} We now consider the configuration (T2).  The analysis will be split in  four further steps.

\textbf{10A.}
To fix ideas, assume that  $\theta(\tau,x)=1$ 
for all $x\in \bigl[b_i(\tau), a_{i+1}(\tau)\bigr]$. 
In this case,  for $x\in [p_k, q_k]$ and $t\in [\tau, \tau+\bar h]$ with $\bar h$ small enough, the semigroup solution 
$U(t,x)$ coincides with the solution to the Cauchy problem with flux $f$, i.e.,
$$U_t +f(U)_x\,=\,0,\qquad U(\tau,x)=u(\tau,x).$$
We shall bound the difference $u-U$, using the fact that $u$
satisfies a conservation law with flux $f$ for all times $t\in K$, and moreover $\tau$ is a 
Lebesgue point in $K$.

Borrowing a technique introduced in \cite{B1}, given $\ve>0$ we partition the interval 
 $\bigl[ p_k, q_k\bigr]$ 
by inserting points
$$b_i(\tau)~=~y_0~<~y_1~<~\cdots~<~y_n~=~ a_{i+1}(\tau)$$
so that 
\bel{tv2}\TV\bigl\{ u(\tau, \cdot)\,;~]y_{j-1}, y_j[\bigr\}~=~u(\tau, y_{j}-)-u(\tau, y_{j-1}+)~<~\ve.\eeq
Note that the above construction is possible because $\theta=1$ and the function $u(\tau,\cdot)$ is increasing.
\v
{\bf 10B.} We now consider the neighborhood around $(\tau,y_j)$. 
Since $\tau\notin{\cal N}$, each $(\tau, y_j)$ is either a point of approximate continuity or a point of approximate jump.
In the neighborhood of $(\tau, y_j)$ we construct an approximation $V_j$, defined as follows.
\begin{itemize} 
\item If $(\tau, y_j)$ is an approximate jump, we define the function
\bel{Vj}V_j(t,x) ~\doteq~\left\{\bega{rl}u^-\quad &\hbox{if}~~x<y_j+\lambda_j (t-\tau),\\
u^+\quad &\hbox{if}~~x>y_j +\lambda_j (t-\tau),\enda\right.\eeq
where 
$$u^\pm\,=\,\lim_{x\to y_j\pm} u(\tau,x),\qquad\qquad \lambda_j\,=\,\frac{ f(u^+)-f(u^-)}{ u^+-u^-}\,.$$
\item 
If  $u(\tau, \cdot)$ is approximately continuous at $y_j$, we simply set
$V_j(t,x)\doteq u(\tau, y_j)$ for all $t,x$.
\bel{Vj2}V_j(t,x) ~\doteq~ u(\tau, y_j) \qquad \mbox{for all} ~ t,x.
\eeq
\end{itemize} 
Then, for any $\lambda^*>0$, we have (see Theorem 2.6 in \cite{Bbook})
\bel{lim8} 
\lim_{t\to \tau+} \frac{1}{ t-\tau} \int_{y_j-\lambda^* (t-\tau)}^{y_j+\lambda^* (t-\tau)} \Big| u(t,x) - V_j(t,x)\Big|\, dx~=~0.
\eeq

{\bf 10C.} We now consider the remaining intervals $[y_{j-1}+ \lambda^* (t-\tau),\, y_j-\lambda^*(t-\tau)]$. 
For each $j$, 
we compare $u(t,x)$ with the 
function 
\bel{Ujflat}U^\flat_j(t,x)~\doteq~u(\tau, x-\lambda_j(t-\tau)\bigr), \qquad\qquad \lambda_j \doteq f' \bigl(u(\tau, y_j-)\bigr). \eeq
Note that $U_j^\flat $ is the solution to the linear Cauchy problem
\bel{Uflin}U_t+ \lambda_j U_x\,=\,0,\qquad\qquad U(\tau,x) \,=\, u(\tau,x).\eeq
Choose $\lambda^*> L_0/\ve$, where $L_0$ is the Lipschitz constant in (\ref{L0}), and consider the domain
\bel{caldj}\D^{(t)}_j~\doteq~\Big\{ (s,x)\,;~~s\in [\tau, t],~~y_{j-1}+ \lambda^*(s-\tau)<x<y_j-
\lambda^*(s-\tau)\Big\}.\eeq
For $(s,x)\in \D^{(t)}_j$ and $s\in K$ we claim that
\bel{uue}
 u(s,x) ~\leq~u(\tau, y_j-)+\ve. \eeq
Indeed, 
by contradiction let us
assume  $u(s,x) >  u(\tau, y_j-) + \ve$.   Then for all $x'\in [x, y_j]$ we have
$$u(s,x') ~\geq ~u(s,x)\geq  u(\tau, y_j-) + \ve~\geq~ u(\tau,x) +\ve.$$
This implies
$$L_0 (s-\tau)~\geq~\int_{y_j-\lambda^*(s-\tau)}^{y_j} \big| u(s,x') - u(\tau,x')\bigr|\, dx' ~\geq~\ve \lambda^*(s-\tau).$$
If we had chosen chose $\lambda^*> L_0/\ve$, this yields  a contradiction.
Similarly, we obtain
$$u(s,x)~\geq~
u(\tau, y_{j-1}+)-\ve~\geq~u(\tau, y_j-) -2\ve. $$
for all $x'\in [y_{j-1} + \lambda^* (s-\tau),\, y_j]$.
Combining the two above estimates, we conclude that, for $s\in K$,
\bel{supu}\sup_{(s,x)\in \D^{(t)}_j}~\bigl| u(s,x)-u(\tau, y_j-)\bigr|~\leq~2\ve.\eeq

We remark that (\ref{supu}) amounts to the {\bf tame oscillation condition} introduced in \cite{BGo}.
In the present setting, this condition is an easy consequence of the local monotonicity of the solution $u(t,\cdot)$,
for $t\in K$.
\v

To estimate the difference $u-U^\flat_j$ 
we  recall a technical  Lemma  (see \cite{Bbook} Lemma 9.3 for a proof).

\begin{lemma}\label{l:31} Let $w:\,]a,b[\,\mapsto\R$ be a bounded function such that, for some positive measure $\mu$, one has
\bel{mubo}
\left| \int_\xi^\zeta w(s)\, dx\right|~\leq~\mu\bigl([\xi,\zeta]\bigr)\qquad\hbox{whenever} ~~a<\xi<\zeta<b.\eeq
Then 
\bel{wbo} \int_a^b \bigl|w(x)\bigr|\, dx~\leq~\mu\bigl(]a,b[\bigr).\eeq
\end{lemma}

\begin{figure}[ht]
\centerline{\hbox{\includegraphics[width=9cm]{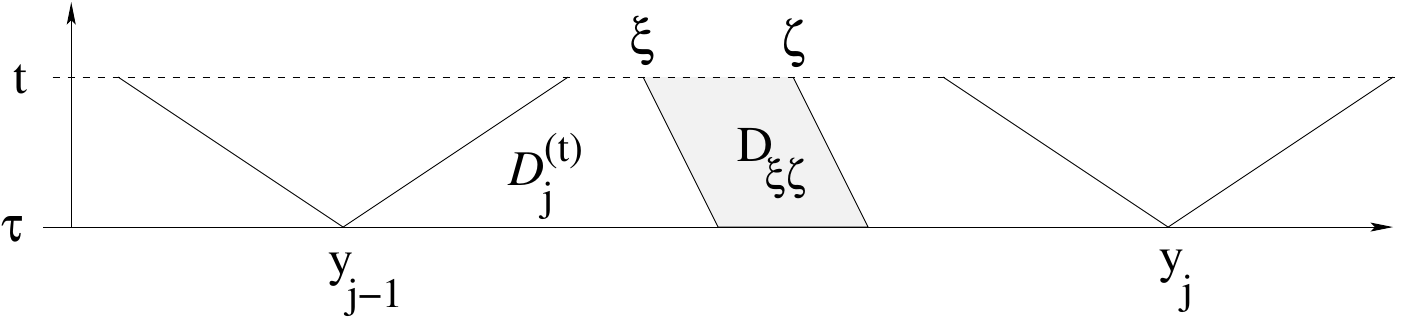}}}
\caption{\small The domain ${\cal D}_j^{(t)} $ at (\ref{caldj}) and the domain $D_{\xi\zeta}$ at (\ref{Dxz}).}
\label{f:df97}
\end{figure}

To apply the lemma, let $\lambda_j$ be as in (\ref{Ujflat}).
Consider any two points $\xi, \zeta$ such that
$$y_{j-1}+ \lambda_j (t-\tau)\, <\, \xi\,<\, \zeta\,<\,y_j- \lambda_j (t-\tau)\,.$$
With reference to Fig.~\ref{f:df97}, we apply the divergence theorem on the domain 
\bel{Dxz}D_{\xi\zeta}~\doteq~\Big\{(s,x)\,;~~s\in [\tau, t],~~\xi+ \lambda_j (s-t)<x<\zeta+ \lambda_j(s-t)\Big\}.\eeq
Note that $D_{\xi\zeta} \subset D_j^{(t)}$.

For notational convenience, in the following estimate we write
$$\xi(s)\,\doteq\,\xi+ \lambda_j (s-t),\qquad\qquad \zeta(s)\,\doteq\,\zeta+ \lambda_j(s-t).$$
Since $u$ satisfies (\ref{1})  and $\theta=1$ when $s\in K$, while $U_j^\flat$ satisfies (\ref{Uflin}),  denoting by
 $F$ the flux function in (\ref{Flux}) we find
\begin{align}
&\int_\xi^\zeta \Big[ u(t,x) - U_j^\flat(t,x)\Big]\, dx
\nonumber \\
 &=
\int_{[\tau,t]\cap K} \Big[ f\bigl( u(s,\xi (s))\bigr)
- f\bigl( u(s,\zeta(s))\bigr) - \lambda_j \bigl( u(s,\xi (s))-u(s,{ \zeta (s)})\bigr)\Big]\, ds
\nonumber \\
 &\quad +\int_{[\tau,t]\setminus K} \Big[F\bigl(s,\xi(s)\bigr)- F\bigl(s,\zeta(s)\bigr)
- \lambda_j \bigl( u(s,\xi (s))-u(s,\zeta (s))\bigr)\Big]\, ds
\nonumber \\
&= \O(1)\cdot\sup_{(s,x)\in \D_j}~\bigl| u(s,x)-u(\tau, y_j-)\bigr| \cdot \int_{[\tau,t]\cap K}  \TV\Big\{ u(s,\cdot)\,;~~\bigl] \xi(s), \zeta(s)\bigr[\, \Big\} \,ds
\nonumber  \\
 &\quad+ 
\O(1) \cdot\int_{[\tau,t]\setminus K} 
\TV\Big\{ F(s,\cdot);\bigl] \xi(s), \zeta(s)\bigr[\, \Big\} +  \TV\Big\{ u(s,\cdot);\bigl] \xi(s), \zeta(s)\bigr[\, \Big\} 
\,ds . 
\label{uUb}
\end{align}

We now introduce two positive measures $\mu_j, \nu_j$ on the interval $$I_j(t)\,\doteq \,\bigl]y_{j-1}+ \lambda^*(t-\tau), 
y_j - \lambda^* (t-\tau)\bigr[\,.$$ 
These are  defined by setting 
\bel{mude}
\mu_j\bigl(]\xi,\zeta[\bigr)~=~ \int_{[\tau,t]\cap K}  \TV\Big\{ u(s,\cdot)\,;~~\bigl] \xi(s), \zeta(s)\bigr[\, \Big\} ds,\eeq
\bel{mu'de}
\nu_j\bigl(]\xi,\zeta[\bigr)=\int_{[\tau,t]\setminus K} 
\left[
\TV\Big\{ F(s,\cdot)\,;~\bigl] \xi(s), \zeta(s)\bigr[\, \Big\} +  \TV\Big\{ u(s,\cdot)\,;~\bigl] \xi(s), \zeta(s)\bigr[\, \Big\} 
\right]
ds\,, \eeq
for every open subinterval $]\xi,\zeta[\subset I_j$.
Using Lemma~\ref{l:31} on every interval $I_j$ and recalling (\ref{supu}), from (\ref{uUb}) we obtain
\bel{uUj} \int_{I_j} \bigl| u(t, x) - U_j^\flat(t,x)\bigr|\, dx~\leq~\O(1)\cdot 
\ve\,\mu_j(I_j) +\O(1)\cdot \nu_j (I_j).\eeq
\v
{\bf 10D.}
Finally, we observe that 
\bel{l30}\sum_j \mu_j(I_j) ~=~\O(1)\cdot (t-\tau) ,\eeq
and
\bel{l31}\sum_j \nu_j(I_j) ~=~\O(1)\cdot \meas\bigl( [\tau,t]\setminus K\bigr),
\qquad\qquad \lim_{t\to \tau+} \frac{1}{ t-\tau}  \sum_j \nu_j(I_j) ~=~0.\eeq
In conclusion,  recalling (\ref{Vj})--(\ref{lim8}), for any $\ve>0$ we obtain
\begin{align}
& \sum_j \int_{I_j} \bigl| u(t, x) - U_j^\flat(t,x)\bigr|\, dx 
+ \sum_j 
 \int_{y_{j}-\lambda^*(t-\tau)} ^{y_j+\lambda^*(t-\tau)} 
\Big|u(t,x) - 
{ V_j} 
(t,x)\Big|\, dx
\nonumber \\
&\qquad\ds =~ \O(1)\cdot \ve (t-\tau) + o(t-\tau)
\label{boo}
\end{align}
Here the last term $o(t-\tau)$  denotes a quantity such that 
$$\lim_{t\to \tau+} \frac{o(t-\tau)}{ t-\tau}~=~0.$$
The same bounds holds replacing $u$ with  the semigroup solution $U(t)= S_{t-\tau} u(\tau)$.
We this conclude
\bel{lim33}\liminf_{t\to \tau+} \frac{1}{ t-\tau} \int_{p_k+\lambda^*(t-\tau)}^{q_k -\lambda^*(t-\tau)} \bigl| u(t,x) - U(t,x)\bigr|\, dx~=~0.
\eeq
 Indeed, this limit is $\leq C\ve$
for some constant $C$ and any $\ve>0$.

\v

\begin{figure}[ht]
\centerline{\hbox{\includegraphics[width=10cm]{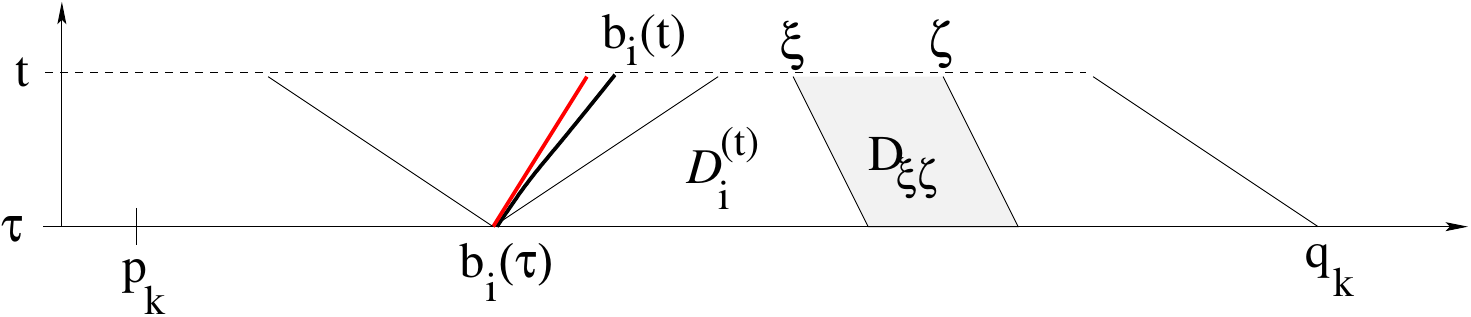}}}
\caption{\small Constructing an approximate solution in case (T3).}
\label{f:df98}
\end{figure}

{\bf 11.} Next, we analyze the behavior of the solution $u$ on an interval $[p_k, q_k]$ of type (T3).  
To fix ideas,
assume that $u(\tau,\cdot)$ is constant for $x<b_i(\tau)$ and decreasing for $x>b_i(\tau)$, 
with a downward jump at $x=b_i(\tau)$.
By assumption, $\bigl(\tau, b_i(\tau)\bigr)$ is a point of approximate jump for $u$, with left and right states
$u^-=u_i(\tau)$, $u^+=u\bigl(\tau, b_i(\tau)+\bigr)$.   Therefore, setting 
$$\lambda\,=\, \frac{g(u^+)- g(u^-)}{ u^+-u^-}$$
the function
\bel{Usharp}
U_i^\sharp(t,x)~\doteq~\left\{\bega{rl} u^-  \qquad &\hbox{if}\quad x< b_i(\tau) + \lambda(t-\tau),\\[1mm]
u^+ \qquad &\hbox{if}\quad x> b_i(\tau) + \lambda (t-\tau),
 \enda\right.\eeq
satisfies
\bel{lim12}  \lim_{t\to\tau+} ~\frac{1}{ t-\tau}\int_{b_i(\tau) - \lambda^*(t-\tau)}^{b_i(\tau) + \lambda^*(t-\tau)} \Big|u(t,x) - U_i^\sharp(t,x)\Big|\, dx~=~0.\eeq
Notice that this implies, in particular, that 
$$b_i(t) ~<~b_i(\tau) + \lambda^*(t-\tau)$$
for all $t\in K$ with $t-\tau$ small enough.
On the domain 
$$D^{(t)}_i ~=~\Big\{ (s,x);~s\in [\tau,t]\cap K,~~b_i(\tau) + \lambda^*(s-\tau)<x<
q_k-\lambda^*(s-\tau)
\Big\}$$
choosing $\lambda^*>L_0/\ve$ we obtain
$$ \bigl| u(t,x) - u^+\bigr| ~ =~
\bigl| u(t,x) -  u(\tau, b_i(\tau)+)\bigr|~\leq~2\ve. 
$$
Setting $\lambda_i\doteq g'(u^+)$ and letting
$U_i^\flat(t,x)~\doteq ~u\bigl(\tau, x -\lambda_i(t-\tau)\bigr)$ be the solution to
$$v_t + g(v)_x=0,\qquad v(\tau,x)= u(\tau,x),$$
the same analysis as in step {\bf 10C} yields the bound (see Figure~\ref{f:df98})
\bel{lim13}  \lim_{t\to\tau+} ~\frac{1}{ t-\tau}\int_{b_i(\tau) + \lambda^*(t-\tau)}^{q_k -\lambda^*(t-\tau)}  \Big|u(t,x) - U_i^\flat(t,x)\Big|\, dx~=~\O(1)\cdot \ve.\eeq
 Finally, over the interval $[p_k, b_i(\tau)- \lambda^*(t-\tau)]$ 
the function $u(t,\cdot)=u_i(t)$ is constant. Moreover,
\bel{omi}
\dot u_i(\tau) = \frac{f(u_i(\tau))-g(u_i(\tau))}{b_i(\tau)-a_i(\tau)}\,.
\eeq
Therefore 
\bel{lim14} \lim_{t\to\tau+} \frac{1}{ t-\tau}\int_{p_k}^{b_i(\tau)- \lambda^*(t-\tau)}\Big|u(t,x) - u_i(\tau) - (t-\tau) \dot u_i(\tau)\Big|\, dx~=~0.\eeq
Since all three limits 
\eqref{lim12}, \eqref{lim13} and \eqref{lim14}
remain valid by replacing the solution $u$ with the semigroup
solution $U(t,\cdot) = S_{t-\tau}u(\tau)$, we conclude that
\bel{lim16}  \lim_{t\to\tau+} \frac{1}{ t-\tau}\int_{p_k}^{q_k- \lambda^*(t-\tau)}\Big|u(t,x) -U(t,x)\Big|\, dx~=~\O(1)\cdot \ve.\eeq
where $\ve$ can be chosen arbitrarily small.

{\bf 12.} Finally, we analyze the behavior of the solution $u$ on an interval $[p_k, q_k]$ of type (T4).  To fix ideas,
we again assume that $u(\tau,\cdot)$ is constant for $x<b_i(\tau)$ and decreasing for $x>b_i(\tau)$.
Given $\ve>0$, we seek a leading order approximation to $u(t,\cdot)$ on the interval $\bigl[p_k+\lambda^*(t-\tau), \,q_k -\lambda^*(t-\tau)\bigr]$, where $\lambda^*= L_0/\ve$.

\begin{figure}[ht]
\centerline{\hbox{\includegraphics[width=9cm]{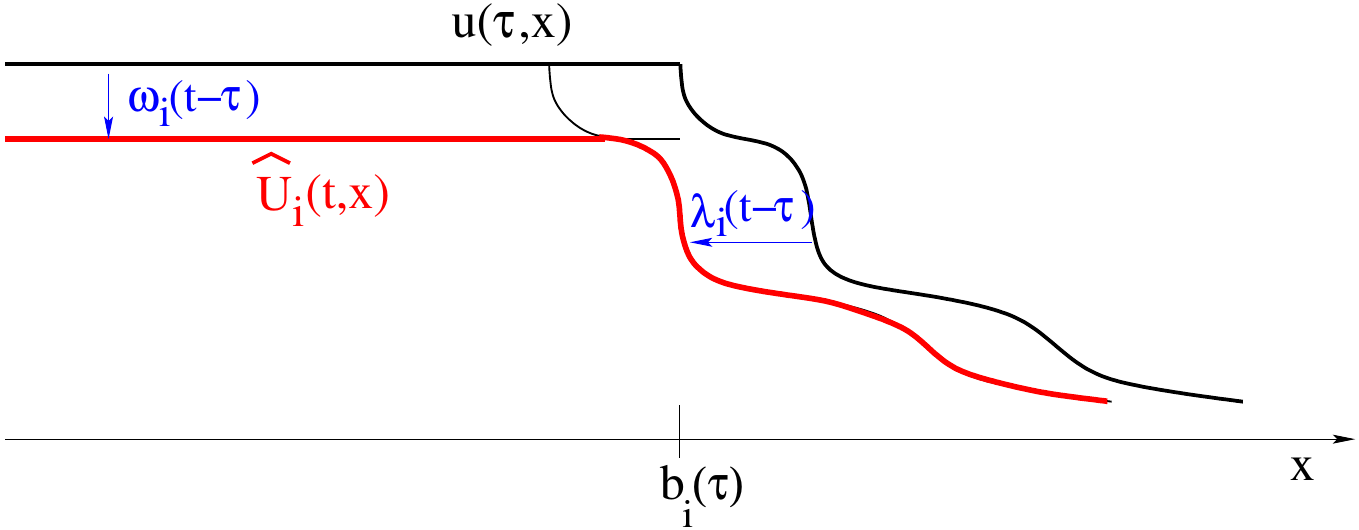}}}
\caption{\small The leading order approximate solution $\Hat U_i$ defined at (\ref{UF}).
Here we are assuming $\lambda_i<0$.
}
\label{f:df100}
\end{figure}

For $t\in K$, setting $\omega_i=\dot u_i(\tau)$ as in (\ref{omi}), we have 
\bel{uupb}u(t,x)~=~u_i(t)~=~u_i(\tau) + (t-\tau)\omega_i  + o(t-\tau)\qquad\qquad \hbox{for} ~~x\leq b_i(t).\eeq
Moreover, $u(t,x)\leq u_i(t)$ for $x\in [p_k, q_k]$.  As in (\ref{uue}) by  the Lipschitz continuity of the map $t\mapsto u(t,\cdot)$ and the fact that $u(t,\cdot) $ is non-increasing in $x$, we obtain
\bel{u11}u(t,x)~\geq~u(\tau, b_i(\tau)) -2\ve\qquad\qquad \hbox{if} \quad t>\tau, \quad p_k<x<q_k-\lambda^*(t-\tau).\eeq
To define a leading order approximation $\Hat U_i(t,\cdot)$ to the solution $u(t,\cdot)$ for $t>\tau$, 
we consider the characteristic speed
\bel{lisp}\lambda_i\,\doteq\, g'\bigl(u(\tau, b_i(\tau)\bigr)\eeq
and define (see Fig.~\ref{f:df100})
\bel{UF}\Hat U_i(t,x)~\doteq~\min\Big\{ u_i(\tau) +\omega_i (t-\tau),~~u(\tau, x-\lambda_i(t-\tau)\bigr)\Big\}.\eeq

\begin{figure}[ht]
\centerline{\hbox{\includegraphics[width=12cm]{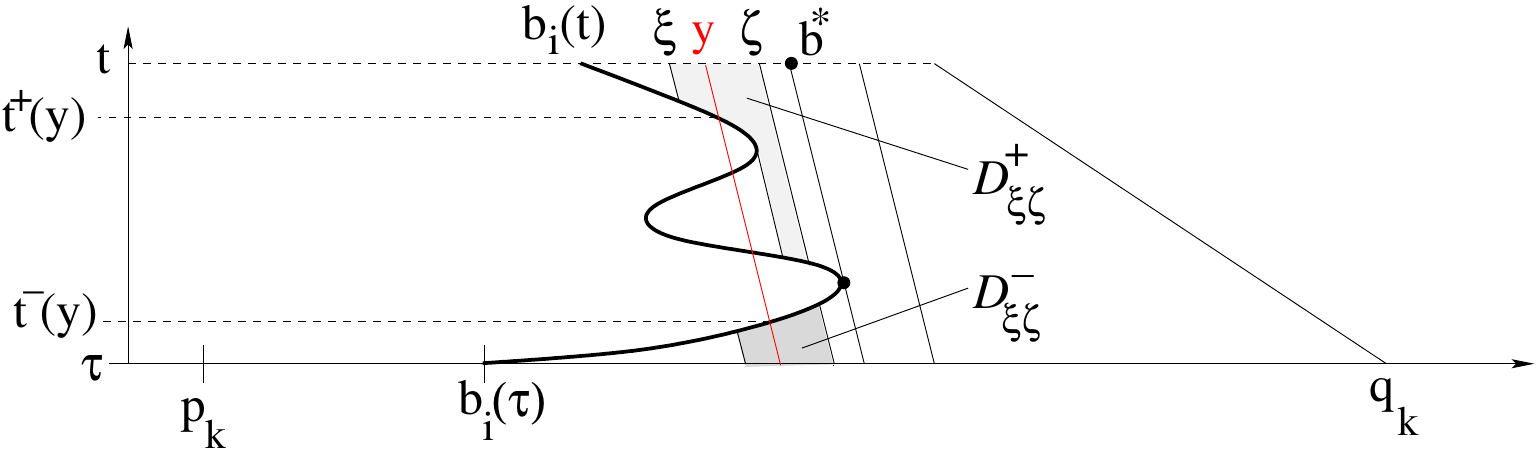}}}
\caption{\small The upper and lower shaded regions are the domains $\D^-_{\xi\zeta}$ and $\D^+_{\xi\zeta}$ in  (\ref{Dxz-}) and (\ref{Dxz+}), used to estimate the error at (\ref{lim15}).}
\label{f:df99}
\end{figure}

We claim that
\bel{lim15}
\lim_{t\to \tau+} \frac{1}{ t-\tau}\int_{p_k}^{q_k - \lambda^*(t-\tau)} \Big|
u(t,x)- \Hat U_i(t,x)\Big|\, dx~=~\O(1)\cdot \ve.\eeq
As shown in Fig.~\ref{f:df99}, consider the point
\bel{Bt}b^*~\doteq~\max_{s\in [\tau,t]\cap K} \Big\{b_i(s)+ \lambda_i(t-s)\Big\}~=~b_i(t^*)+ \lambda_i(t-t^*)\eeq
for some $t^*\in [\tau,t]\cap K$.   Note that at least one such time $t^*$ exists because $K$ is compact and $b_i(\cdot)$ is continuous restricted to $K$.

Trivially, the definition implies
$b^*\geq b_i(t)$.   On the domain
$${\cal D}_i^{(t)} ~\doteq~\Big\{ (s,x)\,;~~s\in [\tau,t]\cap K,~~~x\in \bigl[b^*+\lambda_i(s-t),~q_k-\lambda^*(s-\tau)\bigr]
\Big\}$$
we have $\theta\equiv 0$ and the flux is constantly equal to $g(u)$.
In view of (\ref{u11}),  the same arguments used in step {\bf 10} to analyze the case (T2) now yield
\bel{lim166}
\int_{b^*}^{q_k-\lambda^*(t-\tau)}\bigl|
u(t,x)-\Hat U_i(t,x)\bigr|~=~\O(1)\cdot \ve (t-\tau).\eeq

It remains to study what happens on the interval $[p_k+\lambda^*(t-\tau), b^*]$.
For $y< b^*$ consider the times (see Fig.~\ref{f:df99})
$$t^+(y)~\doteq~\min\Big\{ t'\in [\tau,t]\cap K\,;~~y+\lambda_i(t-s)\geq b_i(s)\quad\forall s\in [t',t]\Big\},$$
$$t^-(y)~\doteq~\max\Big\{ t'\in [\tau,t]\cap K\,;~~y+\lambda_i(t-s)\geq b_i(s)\quad\forall s\in [\tau, t']\Big\}.$$

To complete the proof,  the main idea is that, in the region where $x>b_i(s)$, $s\in K$, the function  $u=u(s,x)$ satisfies the conservation law with flux $g$.    Since the oscillation of $u$ on this domain can be taken very small,  the characteristic speed 
$g'\bigl(u(t,x)\bigr)\approx \lambda_i$ 
is nearly constant.   In turn, the function  $u=u(t,x)$ is nearly constant along lines with slope $\lambda_i$.  
Knowing that along the curve $b_i(\cdot)$  one has 
$u(s, b_i(s))=u_i(s)$, this will yield the desired bounds.

We first apply the divergence theorem on every domain of the form (see the lower shaded region in Fig.~\ref{f:df99})
\bel{Dxz-}D^-_{\xi\zeta}~\doteq~\Big\{ (s,x)\,;~~x= y + \lambda_i(s-t),~~s\in [\tau, t^-(y)],
\quad \hbox{for some}~
y\in [\xi,\zeta]\Big\},\eeq
where $u$ satisfies the conservation law with flux $g$.  
Note that this can be achieved by the integral  formulation (\ref{wsol}), choosing a test function  
$\vp$ which approaches the characteristic function of $\D^-_{\xi\zeta}$.  For example,  for every $y\in [\xi,\zeta]$, 
we can set
$$\vp(s, y+\lambda_i(s-t)\bigr)~=~\left\{ \bega{cl} 0\qquad &\hbox{if} \quad s\notin [0, t^-(y)+\epsilon],\\[1mm]
1\qquad &\hbox{if} \quad s\in [\epsilon , t^-(y)],\\[1mm]
\epsilon^{-1} s\qquad &\hbox{if} \quad s\in [0,\epsilon],\\[1mm]
1-\epsilon^{-1} \bigl(s-t^-(y)\bigr)\qquad &\hbox{if} \quad s\in [t^-(y), t^-(y)+\epsilon].\enda\right.$$
Then we let $\epsilon \to 0+$.

Since the map $t\mapsto u_i(t)$
is strictly decreasing for $t\in K$,
for every $y$ there holds
$$u\Big(t^-(y), y-\lambda_i\bigl(t^-(y)-t\bigr)\Big)~=~u_i\bigl(t^-(y)\bigr) ~\geq~u_i(t) .$$
Proceeding as in step {\bf 10C}, we thus obtain a bound of the form
\bel{bin5}\int_{p_k+\lambda^*(t-\tau)}^{b^*} \Big[ u_i(t)-u(\tau, x-\lambda_i(t-\tau)) \Big]_+\, dx ~=~\O(1)\cdot \ve(t-\tau).\eeq
Here $[z]_+\doteq\max\{z,0\}$ denotes the positive part of a number $z$.
Next, for every $\xi,\zeta$ we apply the divergence theorem on the 
domain
\bel{Dxz+}D^+_{\xi\zeta}~\doteq~\Big\{ (s,x)\,;~~x= y + \lambda_i(s-t),~~s\in [t^+(y),t],
\quad \hbox{for some}~
y\in [\xi,\zeta]\Big\},\eeq
(see the upper shaded region in Fig.~\ref{f:df99}).  Since 
$$u\Big(t^+(y), y-\lambda_i\bigl(t^+(y)-t\bigr)\Big)~=~u_i\bigl(t^+(y)\bigr)~\geq ~u_i(t),$$ 
the same technique now yields
\bel{bin7}\int_{p_k+\lambda^*(t-\tau)}^{b^*} \Big[ u_i(t) - u(t,x)\Big]_+\, dx ~=~\O(1)\cdot \ve(t-\tau).\eeq
By (\ref{bin5}), recalling (\ref{uupb}) and the definition of $\Hat U_i$ at  \eqref{UF}, we obtain
\bel{bin8}\int_{p_k+\lambda^*(t-\tau)}^{b^*} \Big| u_i(\tau)+ (t-\tau)\omega_i-u\bigl(\tau, x-\lambda_i(t-\tau)\bigr) \Big| \, dx ~=~\O(1)\cdot \ve(t-\tau).\eeq
On the other hand, since $u(t,x)\leq u_i(t)$, by (\ref{bin7}) we have
\bel{bin9}\int_{p_k+\lambda^*(t-\tau)}^{b^*} \Big| u_i(\tau)+ (t-\tau)\omega_i-u(t,x) \Big| \, dx ~=~\O(1)\cdot \ve(t-\tau).\eeq
Combining (\ref{bin8}) and (\ref{bin9}), we conclude
\bel{bin77}\int_{p_k+\lambda^*(t-\tau)}^{b^*} \Big| u(t,x)- \Hat U_i(t,x)\Big|\, dx ~=~\O(1)\cdot \ve(t-\tau).\eeq
\v
{\bf 13.} As before, we write $ U(t,\cdot) \doteq  S_{t-\tau} u(\tau,\cdot)$ for the corresponding semigroup
trajectory.
Consider any time $\tau\in K\setminus ({\cal N}\cup {\cal N}')$.   
Summarizing the previous steps, for each $\ve>0$ we can find a finite covering
$[0,1]\subset \bigcup_{k=1}^N ]p_k,q_k[$ and a speed $\lambda^*$ sufficiently large so that
\bel{lim6}\bega{l}\ds  \liminf_{t\to  \tau+} \frac{1}{ t-\tau}  \int_0^1 \Big|u(t,x) -  U(t,x)\Big|\, dx\\[4mm]
\qquad\ds \leq ~\sum_k ~ \lim_{t\to \tau+, \,t\in K_\alpha} ~ \frac{1}{ t-\tau} \int_{p_k+\lambda^*(t-\tau)}^{q_k-\lambda^*(t-\tau)} \Big|u(t,x) - U_k(t,x)\Big|+  \Big|U(t,x) - U_k(t,x)\Big|\, dx\\[3mm]
\qquad \leq~\ve.
\enda\eeq
This yields (\ref{key}), completing the proof.
\endproof
\vs

{\bf Acknowledgment.} The research of A.\,Bressan was partially supported by NSF with
grant  DMS-2306926, ``Regularity and approximation of solutions to conservation laws".


\end{document}